# Integrals related with Rogers Ramanujan continued fraction and q-products.

N. Bagis
Department of Informatics Aristotle University of Thessaloniki 54006 Greece.
bagkis@hotmail.com

M. L. Glasser
Department of Physics Clarkson University Potsdam, NY 13699-58200(UA)
laryg@clarkson.edu

**Abstract.** We present here a way to evaluate a very wide class of integrals relating Ramanujan`s continued fraction and q-product. To do this we explore briefly a differential equation, which relates these two functions.

## 1. Introduction

$$\eta(\tau) := e^{\pi i\tau/12}\prod_{n=1}^{\infty}(1-e^{2\pi i n\tau})$$

denotes the Dedekind eta function.

$${}_2F_1\left[\begin{matrix}a,b\\ c\end{matrix};x\right] := \sum_{k=0}^{\infty}\frac{(a)_k(b)_k}{(c)_k}\frac{x^k}{k!}$$

is Gauss Hypergeometric function.

For $|q|<1$, the “Rogers-Ramanujan continued fraction” is defined as

$$R(q) := \cfrac{q^{1/5}}{1+\cfrac{q}{1+\cfrac{q^2}{1+\cfrac{q^3}{1+\ldots}}}}$$

We also define the q-product

$$f(-q) := \prod_{n=1}^{\infty}\left(1-q^n\right)$$

Also hold the following relations of Ramanujan

$$\frac{1}{R(q)}-1-R(q) = \frac{f(-q^{1/5})}{q^{1/5}f(-q^5)} \tag{1}$$

$$\frac{1}{R^5(q)}-11-R^5(q)=\frac{f^6(-q)}{q\cdot f^6(-q^5)} \tag{2}$$

We will examine the functions *f*, and *R*.

## 2. q-Integrals

**Theorem 1.**

$$\frac{\pi}{3}\int \eta(i\tau)^4 d\tau =$$

$$=R(e^{-2\pi\tau})^{5/6}\sum_{n=0}^{\infty}\frac{(\frac{1}{6})_n^2}{(\frac{7}{6})_n n!}\left(\frac{11+5\sqrt{5}}{2}\right)^n {}_2F_1\left[\begin{matrix}1/6,-n\\5/6-n\end{matrix};\frac{11-5\sqrt{5}}{11+5\sqrt{5}}\right]R(e^{-2\pi\tau})^{5n}$$

**Proof.** From the logarithmic derivative of (see [A2]):

$$R(q)=\frac{\sqrt{5}-1}{2}\exp\left(\frac{1}{5}\int_1^q \frac{f^5(-t)}{f(-t^5)t}dt\right) \tag{3}$$

and the relation (2) we get the differential equation

$$\frac{5R'(q)}{R(q)\left(\frac{1}{R(q)^5}-11-R(q)^5\right)^{1/6}}=f(-q)^4 q^{-5/6}. \tag{4}$$

One can see from (4) that

$$-\int f(-q)^4 q^{-5/6}dq=$$

$$-6R(y)^{5/6}F_1\left[\tfrac{1}{6},\tfrac{1}{6},\tfrac{1}{6},\tfrac{7}{6},\tfrac{1}{2}(11-5\sqrt{5})R(y)^5,\tfrac{1}{2}(11+5\sqrt{5})R(y)^5\right] \tag{5}$$

Where $F_1\left[a,b_1,b_2,c,x,y\right]=\sum_{m=0}^{\infty}\sum_{n=0}^{\infty}\frac{(a)_{m+n}(b_1)_m(b_2)_n}{m!n!(c)_{m+n}}x^m y^n$ is the Appel function. Using this expansion we get the desired result. □

From (5) we also get

$$\int_0^1 f(-q)^4 q^{-5/6}dq=\pi 2^{1/6}(\sqrt{5}-1)^{5/6}{}_2F_1\left[\begin{matrix}1/6,1/6\\1\end{matrix};\frac{1}{2}(-123+55\sqrt{5})\right] \tag{6}$$

$$\int_0^1 f(-q^5)^4 q^{-1/6}dq=\frac{1}{8\cdot 2^{1/6}}\pi(\sqrt{5}-1)^{25/6}{}_2F_1\left[\begin{matrix}5/6,5/6\\1\end{matrix};\frac{1}{2}\left(-123+55\sqrt{5}\right)\right] \tag{7}$$

From (4) and (5) we also get

**Proposition 1.**

$$\int_0^1 \eta(ix)^4 dx = 1/2\left(\frac{\sqrt{5}-1}{2}\right)^{5/6} {}_2F_1\left[\begin{matrix}1/6,1/6\\1\end{matrix};\frac{1}{2}\left(-123+55\sqrt{5}\right)\right] \qquad (8)$$

**Proof.**

setting $a = 1/R(1)^5 = \dfrac{11+5\sqrt{5}}{2}$ and using (4)

$\displaystyle\int_0^1 q^{-5/6} f^4(-q)dq = \int_a^\infty \frac{x^{-5/6}}{(x^2-11x-1)^{1/6}}dx = \int_0^\infty \frac{x^{-1/6}}{(x+a)^{5/6}(x+a+1/a)^{1/6}}dx$, having in mind

that $\displaystyle\int_0^\infty \frac{x^{s-1}}{(x+A)^\mu (x+B)^\nu}dx = \frac{B^{s-\nu}}{A^\mu}B(s,\mu+\nu-s)\,{}_2F_1\left[\begin{matrix}\mu,s\\\mu+\nu\end{matrix};1-\frac{B}{A}\right]$, [GR] we get the result. □

In the same way from the differential equation (4) we get

**Proposition 2.**

$$\int_a^b f(-q^5)^2 f(-q)^2 q^{-1/2} dq =$$

$$= \frac{2i}{5\sqrt{5}}\sqrt{11+2i}\cdot F\left( i\cdot \operatorname{arcsinh}\left(\frac{f(-x^5)^3\sqrt{(11+2i)x}}{f(-x)^3}\right), \frac{117}{125}+\frac{44i}{125}\right)\Bigg|_a^b,$$

where $F$ is the Elliptic$F$ function see [WW].

**Proof.** Let $g(y) = \dfrac{y^{5/3}}{(1-11y^5-y^{10})^{1/3}}$. It is known that

$$f(-q)^4 q^{-5/6} dq = \frac{dR}{R\left(\frac{1}{R^5}-11-R^5\right)^{1/6}}$$

thus $g(R(q))f(-q)^4 q^{-5/6}dq = \dfrac{g(R)dR}{R\left(\frac{1}{R^5}-11-R^5\right)^{1/6}}$

make the change of variable $\left(\dfrac{1}{R^5}-11-R^5\right) = a^6$, after simplification (with Mathematica) we have

$\displaystyle\int \frac{6a^2}{\sqrt{125+22a^6+a^{12}}}da = -\frac{2i}{5\sqrt{5}}\sqrt{11+2i}F\left(i\cdot\operatorname{arccsch}\left(\frac{\sqrt{11+2i}}{a^3}\right), \frac{117}{125}+\frac{44i}{125}\right)$ and

making use of the relation $\dfrac{1}{R(q)^5}-11-R(q)^5 = \dfrac{f(-q)^6}{q\cdot f(-q^5)^6}$, we get the result. □

The above method can be generalized as we show next:

Let $f_k(y)=\left(\dfrac{y^{5/6}}{(1-11y^5-y^{10})^{1/6}}\right)^{6k+5}$, from $f(-q)^4q^{-5/6}dq=\dfrac{dR}{R\left(\dfrac{1}{R^5}-11-R^5\right)^{1/6}}$,

we get $f_k(R(q))f(-q)^4q^{-5/6}dq=\dfrac{f_k(R)dR}{R\left(\dfrac{1}{R^5}-11-R^5\right)^{1/6}}$

set $\left(\dfrac{1}{R^5}-11-R^5\right)=x^6$ then

$$\int_a^b\left(\frac{f(-q^5)}{f(-q)}\right)^{6k+5}f(-q)^4q^kdq=\left(\int_{u(a)^{1/6}}^{u(b)^{1/6}}\frac{6x^{-6k-1}}{\sqrt{125+22x^6+x^{12}}}dx\right)$$

Hence the following theorem holds:

**Theorem 2.**

$$\int_a^b\left(\frac{f(-q^5)}{f(-q)}\right)^{-6k+5}f(-q)^4q^{-k}dq=-\int_{u(a)}^{u(b)}\frac{x^{k-1}}{\sqrt{125+22x+x^2}}dx \qquad (9)$$

where $u(q):=\dfrac{1}{R(q)^5}-11-R(q)^5$.

Next we define the function:

$$F\begin{bmatrix}A,B,C\\ \lambda,\mu,v\end{bmatrix}:=\int_0^\infty\frac{1}{(x+A)^\lambda(x+B)^\mu(x+C)^v}dx \qquad (10)$$

**Observation.**
The values of $F$, for $\lambda=\mu=1/2$ and non negative integer-$v$ are all known functions.

Now set in Theorem 2, $a$ = 0, then $u(a)=\infty$ and (9) becomes

$$\int_0^b\left(\frac{f(-q^5)}{f(-q)}\right)^{6k+5}f(-q)^4q^kdq=\int_{u(b)}^\infty\frac{1}{(x+11+2i)^{1/2}(x+11-2i)^{1/2}x^{k+1}}dx=$$

$$\int_0^\infty\frac{1}{(x+11+2i+u(b))^{1/2}(x+11-2i+u(b))^{1/2}(x+u(b))^{k+1}}dx=$$

$$F\begin{bmatrix}11+2i+u(b),11-2i+u(b),u(b)\\ 1/2,1/2,k\end{bmatrix}=\text{known function of } u(b).$$

Thus we have:

**Theorem 3.**
For $k$ = 0, 1/2, 1, 3/2, 2, … and every 0 < $a$ < 1 the integral

$$\int_0^a \left(\frac{f(-q^5)}{f(-q)}\right)^{6k+5} f(-q)^4 q^k dq,$$

is a known function of $u(q)$. Further we have

$$\int_0^a \left(\frac{f(-q^5)}{f(-q)}\right)^{6k+5} f(-q)^4 q^k dq = F\left[\begin{matrix}11+2i+u(a), 11-2i+u(a), u(a)\\ 1/2, 1/2, k\end{matrix}\right]$$

We can write Theorem 2 as follows

**Theorem 4.**
If $G$ is expandable in power series near 0 and $G(0)=0$ then

$$\int_a^b G\left(\frac{f(-q)^6}{qf(-q^5)^6}\right)\frac{f(-q^5)^5}{f(-q)}dq = -\int_{u(a)}^{u(b)} \frac{G(x)}{x\sqrt{125+22x+x^2}}dx$$

**Examples.**
**1)** Set $k=1$ in Theorem 2, then

$$\int_a^b \frac{f(-q)^5}{f(-q^5)q}dq = -\operatorname{arcsinh}\left(\frac{11+x}{2}\right)\Bigg|_{u(a)}^{u(b)}$$

**2)** Let $B(x) = -\frac{\sqrt{125+22x+x^2}}{125x} - \frac{11\cdot\log(x)}{625\sqrt{5}} +$

$$+\frac{11}{625\sqrt{5}}\log(125+11x+5\sqrt{5}\sqrt{125+22x+x^2})$$

then for $k=-1$, in Theorem 2 we have

$$\int_a^b \frac{f(-q^5)}{f(-q)^7} f(-q)^4 q dq = B(x)\Big|_{u(a)}^{u(b)}$$

**3)** $\int_0^1 f(-q^5)^2 f(-q)^2 q^{-1/2} dq = \frac{2}{5^{3/4}} K\left(\frac{1}{2} - \frac{11}{10\sqrt{5}}\right),$

$K$ is the elliptic integral of the first kind.

**4)** $\int_0^1 f(-q^5)^3 f(-q) q^{-1/3} dq = \frac{2\pi}{11^{2/3}\sqrt{3}} {}_2F_1\left[\begin{matrix}1/3, 5/6\\ 1\end{matrix}; -\frac{4}{121}\right]$

**Theorem 5.**
In general for $0 < k < 1$

$$\int_0^1 \left(\frac{f(-q^5)}{f(-q)}\right)^{-6k+5} f(-q)^4 q^{-k} dq = 11^{k-1}\pi \csc(k\pi)\, {}_2F_1\left[\begin{matrix}\frac{1}{2}-\frac{k}{2}, 1-\frac{k}{2} \\ 1\end{matrix}; -\frac{4}{121}\right]$$

**Proposition 3.**

If $\rho_1$ is the root of $u(\rho_1)=1/2$, then

$$\int_0^{\rho_1} \left(\frac{f(-q^5)}{f(-q)}\right)^{-6k+5} f(-q)^4 q^{-k} dq = \int_{1/2}^{\infty} \frac{x^{k-1}}{\sqrt{125+22x+x^2}} dx$$

and the integral in the right can be expressed in terms of Elliptic functions when $k$ = 0, -1, -2, …. .

**Examples.**

**1)** $$\int_0^{\rho_1} \frac{f(-q^5)^5}{f(-q)} dq = \frac{\log(\frac{1}{2}(7+3\sqrt{5}))}{\sqrt{5}}$$

**2)** $$\int_0^{\rho_1} \frac{f(-q)^{13}}{f(-q^5)^9} q^{3/2} dq =$$

$$\frac{1}{46875}(250\sqrt{1090} + 88i\sqrt{11-2i}E\left[i\cdot \operatorname{arcsinh}(\sqrt{22-4i}), \frac{117}{125}+\frac{44i}{125}\right] -$$

$$(4+66i)\sqrt{11-2i}F\left[i\cdot \operatorname{arcsinh}(\sqrt{22-4i}), \frac{117}{125}+\frac{44i}{125}\right])$$

**3)** Let $\rho_2$ is the root of $u(\rho_2)=1$, then

$$\int_0^{\rho_2} f(-q^5)^2 f(-q)^2 q^{-1/2} dq = -\frac{2i\sqrt{11+2i}}{5\sqrt{5}} F\left[i\cdot \operatorname{arcsinh}\left(\frac{5\sqrt{5}}{\sqrt{11+2i}}\right), \frac{117}{125}+i\frac{44}{125}\right]$$

**Theorem 6.**

If $\rho_3$ is the root of $u(x) = -11+2i$, then for every $k > -1$

$$\int_0^{\rho_3} \left(\frac{f(-q^5)}{f(-q)}\right)^{6k+5} f(-q)^4 q^k dq = \left(\frac{-11-2i}{125}\right)^{1+k} B\left(\frac{1}{2}, 1+k\right) {}_2F_1\left[\begin{matrix}1/2, 1+k \\ 3/2+k\end{matrix}; \frac{117}{125}+\frac{44i}{125}\right],$$

where $B(a,b) = \dfrac{\Gamma(a)\Gamma(b)}{\Gamma(a+b)}$.

**Proof.** Set $A = 11+2i$, $B = 11-2i$, then if $\rho_3$ is the root of $u(x) = -B$ we have $\rho_3 = -0.2302539558379255... - i\cdot 0.1672892791313823...$

and

$$I=\int_0^{\rho_3}\left(\frac{f(-q^5)}{f(-q)}\right)^{-6k+5}f(-q)^4q^{-k}dq=\int_{-B}^{\infty}\frac{x^{k-1}}{\sqrt{125+22x+x^2}}dx=\int_{-B}^{\infty}\frac{x^{k-1}}{\sqrt{(x+A)(x+B)}}dx=$$

$$=\int_0^{\infty}\frac{(y-B)^{k-1}}{y^{1/2}(y-B+A)^{1/2}}dy=\int_0^{\infty}\frac{y^{-1/2}}{(y-B)^{1-k}(y+A-B)^{1/2}}dy=$$

$$\left(\frac{-11-2i}{125}\right)^{1-k}B\left(\frac{1}{2},1-k\right){}_2F_1\left[\begin{matrix}1-k,1/2\\3/2-k\end{matrix};\frac{117+44i}{125}\right]$$, since the next formula holds:

$$\int_0^{\infty}\frac{x^{s-1}}{(x+A)^{\mu}(x+B)^{\nu}}dx=\frac{B^{s-\nu}}{A^{\mu}}B\left(s,\mu+\nu-s\right){}_2F_1\left[\begin{matrix}\mu,s\\\mu+\nu\end{matrix};1-\frac{B}{A}\right]$$, [GR]

**Other Identities.**

**1)** $$\exp\left(5\sqrt{5}\int_0^a\frac{f(-q^5)^5}{f(-q)}dq\right)=\frac{125+11u+5\sqrt{5}\cdot\sqrt{125+22u+u^2}}{(11+5\sqrt{5})u}$$

where $u=u(a)$

**2)** $\int_0^a\frac{f(-q^5)^{11}}{f(-q)^7}qdq=$ known function of $u(a)$

Etc...

One can get

$$\sqrt{5}\int_0^q\frac{f(-t)^5}{f(-t^{1/5})t^{4/5}}dt=-\int_{u(q^{1/5})}^{\infty}\frac{1}{x\sqrt{125+22x+x^2}}dx\Rightarrow$$

$$\exp\left(\sqrt{5}\int_0^q\frac{f(-t)^5}{f(-t^{1/5})t^{4/5}}dt\right)=\frac{125+11w+5\sqrt{5}\sqrt{125+22w+w^2}}{(11+5\sqrt{5})w}\qquad(11)$$

where $w=u(q^{1/5})=\frac{1}{R^5(q^{1/5})}-11-R^5(q^{1/5})$

**Theorem 7.**
Theorem 2 is equivalent with

$$u'(q)=\frac{f(-q)}{f(-q^5)^5}u(q)\sqrt{125+22u(q)+u(q)^2}\qquad(12)$$

Also from the relation

$$f(-q)^4 q^{-5/6} = \frac{R'(q)}{R\left(\frac{1}{R^5}-11-R^5\right)^{1/6}}$$

we get

$$\int_a^b f(-q)^4 G(R(q))dq = \int_{R(a)}^{R(b)} \frac{G(x)}{x\cdot\sqrt[6]{\frac{1}{x^5}-11-x^5}}dx \qquad (13)$$

If we set $G(x)=\left(\frac{1}{x^5}-11-x^5\right)^{1/6}$, we easily get

$$\log(R(q))\Big|_a^b = \int_a^b \frac{f(-q)^5}{f(-q^5)q}dq\,,$$

which is a formula of Ramanujan first proved by G. Andrews [A2]. Note that this relation is expected since we have use it in the beginning to get the differential equation (4)
In the same way as above once can get:

**Theorem 8.**

$$2\int_a^b \left(\frac{f(-q^{1/5})}{q^{1/5}f(-q^5)}\right)^k q^{-4/5}\frac{f(-q)^5}{f(-q^{1/5})}dq = -\int_{y(a)}^{y(b)} \frac{x^{k-1}}{\sqrt{5+2x+x^2}}dx \qquad (14)$$

where

$$y(q) := \frac{1}{R(q)} - 1 - R(q)\,.$$

If we set $k=0$ in (14) we get a proof of the following Ramanujan`s integral first proved by Son [S]:

$$R(q) = \frac{\sqrt{5}-1}{2} - \frac{\sqrt{5}}{1+\frac{3+\sqrt{5}}{2}\exp\left(\frac{1}{\sqrt{5}}\int_0^q \frac{f^5(-t)}{f(-t^{1/5})t^{4/5}}dt\right)}$$

**Theorem 9.**
If $G$ is expandable in power series near 0 and $G(0)=0$ then

$$\int_a^b G\left(\frac{f(-q^{1/5})}{q^{1/5}f(-q^5)}\right) q^{-4/5}\frac{f(-q)^5}{f(-q^{1/5})}dq = -5\int_{y(a)}^{y(b)} \frac{G(x)}{x\sqrt{5+2x+x^2}}dx$$

**Applications**

**1)** If $y(x)=\frac{1}{R(x)}-1-R(x)$, and $x_0$ is the root of $y(x)=0$ then
$x_0 = 0.6816394360211508...$ and

$$1/5\int_{x_0}^{x}\frac{f(-q)^5}{f(-q^5)}q^{-1}dq=\operatorname{arcsinh}\left(\frac{1}{2}\right)-\operatorname{arcsinh}\left(\frac{1+y(x)}{2}\right), \text{ for every } 0<x<1$$

**2)** $$\int_0^x y(q)^k q^{-4/5}\frac{f(-q)^5}{f(-q^{1/5})}dq=5F\begin{bmatrix}1+2i+y(x),1-2i+y(x),y(x)\\1/2,1/2,-k+1\end{bmatrix}$$

As someone can see we open a way to produce new integrals relating Rogers Ramanujan`s continued fraction with $f$.

## References

**[A1]:** George E. Andrews "Number Theory". Dover Publications, Inc. New York. p.370

**[A2]:** G.E. Andrews, Amer. Math. Monthly, 86, 89-108 (1979)

**[Ap]:** Tom Apostol. "Introduction to Analytic Number Theory". Springer-Verlag, New York, Berlin, Heidelberg, Tokyo 1976, 1984.

**[B1]:** B.C. Berndt. "Ramanujan`s Notebooks Part I". Springer-Verlag, New York Inc.1985.

**[B2]:** B.C. Berndt. "Ramanujan`s Notebooks Part II". Springer-Verlag, New York Inc.1989.

**[B3]:** B.C. Berndt. "Ramanujan`s Notebooks Part III". Springer-Verlag, New York Inc.1991

**[GR]:** I.S. Gradshteyn and I.M. Ryzhik, "Table of Integrals, Series and Products". [Academic Press (1980); (3.197 (1))

**[S]:** S. H. Son. "Some integrals of theta functions in Ramanujan`s lost notebook, Proc. Canad. No. Thy Assoc. No. 5 (R. Gupta and K. S. Williams, eds.), Amer. Math. Soc., Providence.

**[T]:** E. C. Titchmarsh, "Introduction to the theory of Fourier integrals". Oxford University Press, Amen House, London, 1948.

**[WW]:** E.T. Whittaker and G.N. Watson, A Course on Modern Analysis, [Cambridge U.P. 1927].

**[Z]:** I.J. Zucker, SIAM J. Math. Anal. 10.192 (1979)